\title{On the gap between representability and collapsibility}

\documentclass[11pt]{article}
\usepackage{a4}
\usepackage{amssymb,amsmath,amsthm}
\usepackage{graphicx}

\newtheorem{theorem}{Theorem}[section]
\newtheorem{lemma}[theorem]{Lemma}

\newtheorem{proposition}[theorem]{Proposition}

\newcommand{\heading}[1]{\medskip\par\noindent{\bf #1}}
\newcommand{\cmt}[1]{\ifhmode\newline\fi{\sf *** \ \ #1 \\}}

\makeatletter
\newcommand{\ProofEndBox}{{\ifhmode\unskip\nobreak\hfil\penalty50 \else
          \leavevmode\fi\quad\vadjust{}\nobreak\hfill$\Box$
            \finalhyphendemerits=0 \par}}%
\makeatother

\newcommand{\proofend}{\ProofEndBox\smallskip}

\newcommand{\set}[1]{\ensuremath{\left\{#1 \right\}}}
\newcommand{\famf}[1]{\mathcal{#1}}

\newcommand{\scxf}[1]{{\sf #1}}
\newcommand{\nerv}[1]{\scxf{N}(#1)}
\newcommand{\collaps}[1]{\searrow \kern-3mm {}^{#1} \kern2mm}
\newcommand{\cols}[1]{\cdim_0(#1)}

\newcommand\ldim{\lambda} 
\newcommand\cdim{\gamma}  
\newcommand\rdim{\rho} 

\newcommand{\Hred}{\tilde H}
\newcommand{\HrFree}[2]{\Hred_{#1} ({#2})}

\newcommand{\conv}{\mathop{\mathrm{conv}}}
\newcommand{\er}{\mathbb{R}}

\newcommand{\kve}{\mathbb{Q}}

\newcommand\mes{{\mathrm{mes}}}

\newcommand{\R}{\er}

\newcommand\FF{{\mathcal{F}}}
\newcommand\GG{{\mathcal{G}}}
\newcommand\HH{{\mathcal{H}}}
\newcommand\MM{{\mathcal{M}}}
\newcommand\NN{{\mathcal{N}}}

\newcommand\bM{{\bf M}}
\newcommand\bK{{\bf K}}

\newcommand\K{{\scxf K}}
\renewcommand\L{{\scxf L}}

\begin{document}
\author{
{\sc Ji\v{r}\'i Matou\v{s}ek\ \ \ \ \ \ \ \ \ \ \ 
Martin Tancer}\\
   {\footnotesize Department of Applied Mathematics and}\\[-1.5mm]
   {\footnotesize Institute of Theoretical Computer Science (ITI)}\\[-1.5mm]
   {\footnotesize  Charles University, Malostransk\'{e} n\'{a}m. 25}\\[-1.5mm]
{\footnotesize  118~00~~Praha~1,
   Czech Republic}
}
\maketitle

\begin{abstract}
A simplicial complex $\K$ is called \emph{$d$-representable}
if it is the nerve of a collection of convex sets in~$\er^d$;
$\K$ is \emph{$d$-collapsible} if it can be reduced to
an empty complex by repeatedly removing a face
of dimension at most $d-1$ that is contained in a unique maximal
face; and $\K$ is \emph{$d$-Leray} if
every induced subcomplex of $\K$ has vanishing 
homology of dimension $d$ and larger.

It is known that $d$-representable implies $d$-collapsible
implies $d$-Leray, and no two of these notions coincide
for $d\ge 2$.
The famous  Helly theorem 
and other important results in discrete geometry 
can be regarded as results about
$d$-representable complexes, and
in many of these results ``$d$-representable''
in the assumption can be replaced by ``$d$-collapsible''
or even ``$d$-Leray''.

We investigate ``dimension gaps'' among these notions,
and we construct, for all $d\ge 1$, a
$2d$-Leray complex that is not $(3d-1)$-collapsible
and a $d$-collapsible complex that is not
$(2d-2)$-representable.  In the proofs we
obtain two results of independent interest:
(i) The nerve of every finite family of
sets, each of size at most $d$,
is $d$-collapsible. (ii) If the nerve of a simplicial
complex $\K$ is $d$-representable, then $\K$ embeds
in $\er^d$.
\end{abstract}

\section{Introduction}

\heading{\boldmath $d$-representability. }
Helly's theorem \cite{h-umkkm-23}
asserts that if $C_1,C_2,\ldots,C_n$
are convex sets in $\R^d$, $n\ge d+1$ and every $d+1$ 
of the $C_i$ have a common point, then 
$\bigcap_{i=1}^n C_i\ne \emptyset$. This famous theorem and
many others in discrete geometry
deal with
\emph{intersection patterns} of convex sets in $\R^d$,
and they can be restated using the notion of
\emph{$d$-representable} simplicial complexes.

We recall that the \emph{nerve} $\nerv{\famf{S}}$ of a family 
$\famf{S} = \set {S_1, S_2, \dots, S_n}$ is the
simplicial complex with vertex set $[n]:=\{1,2,\ldots,n\}$
and with a set $\sigma\subseteq [n]$ forming a simplex
if $\bigcap_{i\in \sigma} S_i\ne\emptyset$.
A simplicial complex $\K$ is
\emph{$d$-representable} if it is 
isomorphic to the nerve of a family of convex sets in $\er^d$
(all simplicial complexes throughout this paper
are assumed to be finite). 

In this language Helly's theorem
implies that a $d$-representable complex is determined
by its $d$-skeleton.
Other examples of theorems that can be seen as statements
about $d$-representable complexes include 
the fractional Helly theorem of Katchalski and Liu 
\cite{kl-pgr-79}, the colorful Helly theorem of Lov\'asz
(\cite{Lovasz-colhelly}; also see \cite{b-gct-82}),
the $(p,q)$-theorem of Alon and Kleitman \cite{ak-pcs-92a},
and the Helly-type result of Amenta \cite{a-spiht-96}
(conjecture by Gr\"unbaum and Motzkin). 
Among the deepest results concerning $d$-representable
complexes is a complete characterization of their
$f$-vectors\footnote{The \emph{$f$-vector} of
a $d$-dimensional simplicial complex $\K$ is the
integer vector $(f_0,f_1,\ldots,f_d)$,
where $f_i$ is the number of $i$-dimensional simplices in $\K$.}
conjectured by Eckhoff and proved by
Kalai \cite{Kalai-suffeckhoff,Kalai-nececkhoff}.
We also refer to
\cite{DanzerGrunbaumKlee,e-hrctt-93,Mat-dg} for more
examples and background.

\heading{\boldmath $d$-collapsibility and $d$-Leray complexes. }
Wegner in his seminal 1975 paper \cite{w-dcnfc-75} 
introduced $d$-collapsible simplicial complexes.
To define this notion, we first introduce an
\emph{elementary $d$-collapse}. Let
 $\K$ be a simplicial complex and let $\sigma,\tau\in\K$ be
faces (simplices) such that

(i) $\dim\sigma\le d-1$,

(ii) $\tau$ is an inclusion-maximal face of $\K$,

(iii) $\sigma\subseteq\tau$, and

(iv) $\tau$ is the \emph{only} face of $\K$
satisfying (ii) and (iii).

\begin{figure}
\centering
\mbox{
\includegraphics[width=\textwidth]{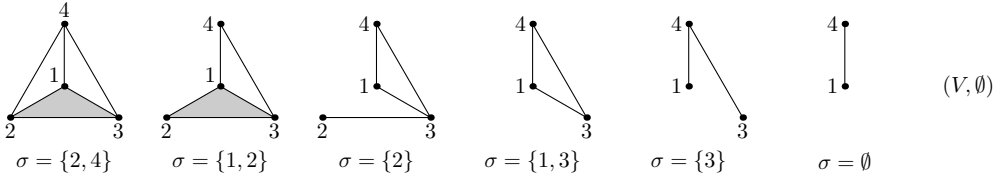}
}
\caption{An example of $2$-collapsing.}
\label{FigCol}
\end{figure}

\noindent
Then we say that $\sigma$ is a \emph{$d$-collapsible face}
of $\K$ and that the simplicial complex
$\K':=\K\setminus\{\eta\in\K: \sigma \subseteq\eta \subseteq \tau\}$
arises from $\K$ by an elementary $d$-collapse.
A~simplicial complex $\K$ is \emph{$d$-collapsible}
 if there exists a sequence of elementary $d$-collapses
that reduces $\K$ to the empty complex $\emptyset$.
Fig.~\ref{FigCol} shows an example of $2$-collapsing.

Another related notion is a \emph{$d$-Leray}
simplicial complex, where $\K$ is $d$-Leray if
every induced subcomplex of $\K$ (i.e., a subcomplex
of the form $\K[X]:=\{\sigma\cap X: \sigma\in\K\}$
for some subset $X$ of the vertex set $V(\K)$)
has zero homology (over $\kve$) in dimension $d$ and larger.

Wegner \cite{w-dcnfc-75} proved that
every $d$-representable complex 
is $d$-collapsible and every $d$-collapsible
complex is $d$-Leray. 
By inspecting  proofs of several theorems
about intersection patterns of convex sets in $\R^d$,
i.e. about $d$-representable complexes, one 
can sometimes see that they actually use only $d$-collapsibility,
and thus they are valid for all $d$-collapsible complexes
(good examples, among those mentioned earlier, are
the fractional Helly theorem
and the colorful Helly theorem).

With more work it has been shown that all of the
results mentioned above and some others also hold
for $d$-Leray complexes. For example, for Helly's theorem
this follows essentially from Helly's own
topological generalization \cite{h-usvam-30},
for the $(p,q)$-theorem this was proved in \cite{AKMM01},
and for the colorful Helly theorem
and for Amenta's theorem this was shown recently
by Kalai and Meshulam \cite{KalMes1,KalMes2}.
Kalai's characterization of $f$-vectors of $d$-representable
complexes is also valid for the $f$-vectors
of $d$-Leray complexes, showing that $f$-vectors
cannot distinguish these classes.

\begin{figure}
\centering
\includegraphics{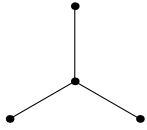}
\caption{A complex that is $1$-collapsible but not $1$-representable.}
\label{Fig3st}
\end{figure}

\begin{figure}
\centering
\includegraphics{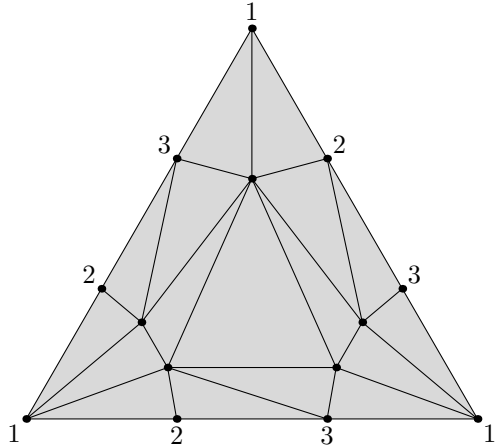}
\caption{The triangulation of the dunce hat~\cite{w-dcnfc-75}; 
vertices with the same numbers should be identified.}
\label{FigDH}
\end{figure}

These results indicate that the notions of 
$d$-representable, $d$-collapsible, and $d$-Leray
are similar in some important respects.
However, no two of them coincide. Fig.~\ref{Fig3st} shows
an example of a $1$-collapsible complex
that is not $1$-representable. Wegner \cite{w-dcnfc-75} 
noted that well-known examples 
of $2$-dimensional complexes that are 
contractible but not collapsible, such as suitable
triangulations of  the ``dunce hat'' (Fig.~\ref{FigDH})
or Bing's house (see, e.g., \cite{Hatcher}),
are $2$-Leray but not $2$-collapsible.

\heading{Results. }
The goal of the present paper is to exhibit stronger
differences among these notions; more precisely, to
investigate ``dimension gaps''. We set
$$
\begin{array}{rcll}
\rdim(\K)&:=& \min\{d: \K\mbox{ is $d$-representable}\}
& \mbox{(``representability'')},\\
\cdim(\K)&:=& \min\{d: \K\mbox{ is $d$-collapsible}\}
& \mbox{(``collapsibility'')},\\
\ldim(\K)&:=& \min\{d: \K\mbox{ is $d$-Leray}\}
& \mbox{(``Leray number'')}.\\
\end{array}
$$

\begin{theorem}\label{t:}\ 

{\rm (a) } For every $d\ge 1$ there exists a complex $\K$
with $\cdim(\K)=d$ and $\rdim(\K)=2d-1$ (i.e.,
$d$-collapsible and not $(2d-2)$-representable).

{\rm (b) } For every $d\ge 1$ there exists a complex $\K$
with $\ldim(\K)=2d$ and $\cdim(\K)=3d$ (i.e.,
$2d$-Leray and not $(3d-1)$-collapsible).
\end{theorem}

In part (a), our example is the \emph{nerve} of
a $d$-dimensional  simplicial complex $\L$
that is not embeddable in $\R^{2d-2}$.
A well known example of such $\L$ is the $d$-skeleton
of the $(2d+2)$-dimensional simplex, due to
Van Kampen \cite{k-ker-32} and Flores \cite{f-undkd-}.
The proof of Theorem~\ref{t:}(a) then follows immediately
from the two propositions below, which may be of independent
interest.

\begin{proposition}\label{p:repr-emb}
Let $\L$ be a simplicial complex such that the nerve $\nerv{\L}$ 
is $d$-representable. Then $\L$ embeds in $\er^d$, even
linearly.
\end{proposition}

\begin{proposition}\label{PropCol}
Let $\famf F$ be a finite family of sets, 
each of size at most $d$. Then the
nerve $\nerv{\famf F}$ is $d$-collapsible.
\end{proposition}

For part (b) of Theorem~\ref{t:}, our example is
a $d$-fold join 
the dunce hat triangulation from Fig.~\ref{FigDH}.

\heading{Open problems. } The main questions, which
we unfortunately haven't solved, are:
Can representability be bounded in terms of collapsibility
(formally, is there a function $f_1$ such that
$\rdim(\K)\le f_1(\cdim(\K))$ for all $\K$)?
Can  collapsibility be bounded in terms of
the Leray number
(formally, is there a function $f_2$ such that
$\cdim(\K)\le f_2(\ldim(\K))$ for all $\K$)?
Theorem~\ref{t:} shows $f_1(d)\ge 2d-1$ and
$f_2(2d)\ge 3d-1$.

It is clear that our method cannot give a better lower bound
for $f_1$ than Theorem~\ref{t:}(a), since every $d$-dimensional
complex embeds in $\R^{2d-1}$. A $2$-collapsible complex
whose representability might perhaps be unbounded
was noted by Alon et al. \cite{AKMM01}, namely,
a finite projective plane (regarded as a simplicial
complex, where the lines of the 
projective plane are the 
maximal simplices). More generally,
any \emph{almost-disjoint} set system
is easily seen to be $2$-collapsible, and it would be interesting
to decide whether all  almost-disjoint systems
are $d_0$-representable
for some constant $d_0$.

\section{Representability of the Nerve and Embeddability}

In this section we will prove Proposition~\ref{p:repr-emb}.
First we recall a classical lemma of Radon (\cite{r-mkkde-21};
also see, e.g., \cite{e-hrctt-93}
or \cite{Mat-dg}) in the following form:

\begin{lemma} \label{ThmRad}
Let $P$ be a set of affinely dependent points in $\er^d$. 
Then there exist two 
disjoint affinely independent subsets $A,B\subset P$ with
 $\conv(A) \cap \conv(B) \neq \emptyset.$
\end{lemma}

We will also need the following result of a similar flavor:

\begin{lemma} \label{LemRad}
Let $A$ and $B$ be finite subsets of $\er^d$. Suppose that there is a point $x \in (\conv(A) \cap \conv(B))\setminus \conv(A \cap B)$. 
Then there exist
disjoint affinely independent
sets $A'\subseteq A$ and~$B'\subseteq B$ such that
$\conv (A') \cap \conv (B') \neq \emptyset$.
\end{lemma}

\begin{proof} The proof is similar to the usual
proof of Radon's lemma, only slightly more complicated.

We can write $x$ as a convex combination of points of $A$:
\begin{equation}\label{e:x1}
x=\sum_{a\in A} \alpha_a a,
\end{equation} 
where $\alpha_a\ge 0$ for all $a\in A$ and 
$\sum_{a\in A}\alpha_a=1$.
Similarly 
\begin{equation}\label{e:x2}
x=\sum_{b\in B} \beta_b b
\end{equation}
where $\beta_b\ge 0$ for all $b\in B$ and  
$\sum_{b\in B}\beta_b=1$.
Let $K := A \cap B$ and let 
$$
K^+ := \{p\in K: \alpha_p >\beta_p\},\ \ \ K^-:= K\setminus K^+.
$$
We define the sets 
$A_0:= A\setminus K^-$ and $B_0:= B\setminus K^+$,
and we note that $A_0\cap B_0=\emptyset$ and 
$A_0\cup B_0=A\cup B$.
We claim that $\conv(A_0)\cap\conv(B_0)\ne\emptyset$;
this will imply the lemma, since the desired
affinely independent $A'$ and $B'$ 
can be obtained from $A_0$ and $B_0$ by removing 
redundant points.

For notational convenience we
extend the definition of $\alpha_p$ and $\beta_p$
to all $p\in A\cup B$ by letting $\alpha_p=0$ for $p\not\in A$
and $\beta_p=0$ for $p\not\in B$.
By subtracting (\ref{e:x2}) from (\ref{e:x1}) and rearranging
we get
$$
\sum_{p\in A_0} (\alpha_p -\beta_p)p =
\sum_{p\in B_0} (\beta_p -\alpha_p)p.
$$
All coefficients  on both sides of this equation are nonnegative.
Let us set $S:= \sum_{p\in A_0}(\alpha_p -\beta_p)$. 
Since $\sum_{p\in A}\alpha_p=\sum_{p\in B}\beta_p=1$,
 we also have $S=\sum_{p\in B_0}(\beta_p -\alpha_p)$.
Moreover, since $x\not\in\conv(K)$, at least one
$\alpha_p$ with $p\in A\setminus K$ is nonzero, and thus
$S\ne 0$. We set
$$
y:=\frac1S \sum_{p\in A_0} (\alpha_p -\beta_p)p=
    \frac1S \sum_{p\in B_0} (\beta_p -\alpha_p)p;
$$
thus, $y$ is expressed as a convex combination of
points of $A_0$ and also as a convex combination of
points of $B_0$. Hence $\conv(A_0)\cap\conv(B_0)\ne\emptyset$
as claimed.
\end{proof}

\heading{Proof of Proposition~\ref{p:repr-emb}. }
Let $\L$ be a simplicial complex such that $\nerv{\L}$
is $d$-representable. This means that there
exists a system $(C_\sigma:\sigma\in \L)$ of convex
sets in $\R^d$ such that for every 
collection $\MM\subseteq \L$ of simplices we have 
$\bigcap_{\sigma\in\MM} C_\sigma=\emptyset$ iff
$\bigcap\MM=\emptyset$.

For every $v \in V(\L)$ 
we fix a point  $p(v)\in\bigcap_{\tau\in \L: v \in \tau}
C_{\tau}$ (this intersection is nonempty since
$v\in \bigcap \{\tau\in\L: v \in \tau\}$).

This defines a mapping $p\colon V(\L)\to\R^d$.
For every $\sigma\in\L$ we set
$$
D_\sigma := \conv(p(\sigma)).
$$
We claim that each $D_\sigma$ is a simplex in $\R^d$
and that the $D_\sigma$ form a geometric representation
of $\K$ in $\R^d$. 
To this end, it suffices to verify that the set $p(\sigma)$
is affinely independent for every $\sigma\in \L$,
and that $D_\sigma\cap D_\tau=D_{\sigma\cap\tau}$
for every two simplices $\sigma,\tau\in\L$.

First, let us suppose for contradiction that
$p(\sigma)$ is affinely dependent for some 
$\sigma\in \L$. Then by Radon's lemma
(Lemma~\ref{ThmRad}) there are two disjoint affinely independent 
subsets $A,B\subset p(\sigma)$ with intersecting convex hulls.
Then we have $A=p(\alpha)$ and $B=p(\beta)$ for
disjoint simplices $\alpha,\beta\in \L$.
But we have $p(v)\in C_\alpha$
for all $v\in \alpha$, hence 
$D_\alpha =\conv(p(\alpha))\subseteq C_\alpha$,
and similarly $D_\beta\subseteq C_\beta$.
Then $C_\alpha\cap C_\beta\supseteq D_\alpha \cap
D_\beta\ne\emptyset$, and this contradicts the
assumption that the $C_\sigma$ form a representation
of $\nerv{\L}$. So each $p(\sigma)$ is affinely independent.

Next, let $\sigma, \tau \in\L$. We clearly
have $D_{\sigma \cap \tau} \subseteq D_\sigma \cap D_\tau$.
To prove the reverse inclusion, we assume for contradiction
that there is some $x \in (D_\sigma \cap
D_\tau)\setminus D_{\sigma \cap \tau}$.
Lemma~\ref{LemRad} provides disjoint 
$\sigma' \subseteq \sigma$ and $\tau' \subseteq \tau$ 
with  $D_{\sigma'} \cap D_{\tau'} \neq \emptyset$,
and this is a contradiction as above. 
\proofend

\section{\boldmath $d$-Collapsibility of the Nerve}

Here we prove Proposition~\ref{PropCol}.

Let us assume that the ground set of $\FF$ is $[n]$.
Let us fix an arbitrary linear ordering $\le$ on $\FF$.
The nerve $\bK=\nerv{\FF}$ consists of all intersecting subfamilies
of $\FF$. For $i=1,2,\ldots,n$ let $\bK_i$ consist of all
intersecting families $\GG\in\bK$ with $\min\bigcap\GG=i$
(so the $\bK_i$ form a partition of $\bK$).

Let us consider a $\GG\in \bK_i$. Each of the elements
$1,2,\ldots,i-1$ is excluded from $\bigcap \GG$ by at least
one $G\in\GG$. Let us define the \emph{minimal exclusion
sequence} $\mes(\GG)=(G_1,G_2,\ldots,G_{i-1})$
as follows. First we choose $G_1$ as the smallest set
of $\GG$ with $1\not\in G$. Having already defined sets
$G_1,\ldots,G_{j-1}\in\GG$
(not necessarily all distinct), we define $G_j$ as follows:
If at least one of the sets among $G_1,\ldots,G_{j-1}$ avoids
the element $j$, we let $G_j$ be such a $G_k$ with the smallest possible $k$.
In this case we call $G_j$ \emph{old at $j$}.
On the other hand, if all of $G_1,\ldots,G_{j-1}$ contain $j$,
then we let $G_j$ be the smallest set of $\GG$ not containing $j$,
and we call it \emph{new at $j$}.

Let $M(\GG)\in\bK_i$ be the family consisting of all sets
$G_j$ that occur in $\mes(\GG)$. In particular,
for $i=1$ we have $M(\GG)=\emptyset$ for all $\GG\in\bK_1$.
It is easily seen that we always have $|M(\GG)|\le d$ (indeed,
$G_1$ covers at most $d-1$ elements among $1,2,\ldots,i-1$,
and only these elements may contribute $G_j$'s distinct from $G_1$).
We also note that $\mes(M(\GG))=\mes(\GG)$.

We let $\bM_i=\{M(\GG): \GG\in\bK_i\}$ and $\bM=\bigcup_{i=1}^n \bM_i$.
The families in $\bM$ will be the $d$-collapsible
simplices we will use for $d$-collapsing
the simplicial complex $\bK=\nerv{\FF}$.
We order them first by \emph{decreasing} $i$,
i.e., $\bM_n$ comes first, then $\bM_{n-1}$, etc.,
and within each $\bM_i$ we order the
families lexicographically by their minimal exclusion sequences.
Let $\preceq$ denote this linear ordering of $\bM$.
This defines the sequence of elementary collapses.

Clearly, each simplex $\GG\in\bK$ contains at least
one simplex of $\bM$, namely, $M(\GG)$.
It remains to verify that each $\MM\in \bM$ is contained
in a unique maximal simplex in the simplicial complex
obtained from $\bK$ by collapsing all $\NN \prec \MM$.

We inductively define
$$
\bK_\MM =
\left\{\HH\in \bK \setminus \bigcup_{\NN\prec \MM}\bK_\NN:
\MM\subseteq \HH\right\};
$$
as we will see, this is the set of all simplices removed
from the current simplicial complex by collapsing $\MM$.

We need to express $\bK_\MM$ as the set of all
simplices of $\bK \setminus \bigcup_{\NN\prec \MM}\bK_\NN$
that contain $\MM$ and are contained in a suitable
maximal simplex $T(\MM)$.
As we will see, the desired $T(\MM)$ can be described
as follows (here $\MM\in \bK_i$ and
$\mes(\MM)=(G_1,\ldots,G_{i-1})$):
$$
T(\MM)=\MM\cup \Bigl\{F\in \FF: i\in F\mbox{ and }F>G_j
\mbox{ for all $j\not\in F$ such that $G_j$ is new at $j$}\Bigr\}.
$$
That is, $T(\MM)$ consists of $\MM$ plus those sets of $\FF$ that contain $i$
and satisfy $\mes(\MM\cup\{F\})=\mes(\MM)$.
Clearly $T(\MM)\in\bK$ and $\MM\subseteq T(\MM)$.
We also have $M(T(\MM))=\MM$.

We let $\bK'_\MM=\{\HH\in \bK\setminus\bigcup_{\NN\prec \MM}\bK_\NN
: \MM\subseteq \HH\subseteq T(\MM)\}$,
and by induction we prove
that $\bK'_\MM=\bK_\MM$. This will show that $T(\MM)$ is indeed
the unique maximal simplex containing $\MM$.
So we consider some $\MM\in \bM_i$ and we assume
that $\bK'_\NN=\bK_\NN$
for all $\NN\prec\MM$.

By just comparing the definitions we immediately obtain
$\bK'_\MM\subseteq \bK_\MM$. For the reverse inclusion
let us consider an $\HH\in\bK_\MM$, and for contradiction
let us suppose that $\HH$ contains a set $F\not\in T(\MM)$.
We will exhibit an $\NN\in\bM$ with $\NN\subseteq\HH$
and $\NN\prec\MM$; this will lead to the desired
contradiction, since then we either have
$\HH\in\bK_\NN$ or $\HH$ has been collapsed even earlier.

By the definition of $T(\MM)$, there can be two reasons for
$F\not\in T(\MM)$. First, it might happen that $i\not\in F$.
But then $\min\bigcap\HH>i$, and therefore $\NN:=M(\HH)\prec \MM$.
This leads to a contradiction as explained above, and
hence we may assume $i\in F$.

Second, letting $\mes(\MM)=(G_1,G_2,\ldots,G_{i-1})$,
there may be some $j<i$, $j\not\in F$ such that $G_j$ is
new at $j$ and $F<G_j$ (we cannot have $F=G_j$ since
we assumed $F\not\in\MM$). Let $j$ be the smallest possible
with this property. We consider the family
$\GG=\{G_1,G_2,\ldots,G_{j-1},F,G_{j+1},\ldots,G_{i-1}\}$,
which is in $\bK_i$. Hence $\NN=M(\GG)$ is in $\bM_i$,
and we also have $\NN\subseteq\HH$. It now suffices
to verify that
$\NN\prec \MM$. To this end, we check that the first $j$ terms
of the sequence $\mes(\NN)$ are $G_1,G_2,\ldots,G_{j-1},F$,
since then $\mes(\NN)$ is indeed lexicographically smaller than
$\mes(\MM)=(G_1,G_2,\ldots,G_{j-1},G_j,\ldots)$.
   
Let us suppose that $\mes(\NN)$ agrees with
$(G_1,G_2,\ldots,G_{j-1},F)$ in the first $k-1$ terms,
and we want to check that the $k$th terms agree as well.
This is clear if $k<j$ and  $G_k$ is old at $k$
in $\mes(\MM)$. If $k<j$ and $G_k$ is new at $k$
in $\mes(\MM)$, then $k\in F$ or $F>G_k$, for otherwise,
we should have taken $k$ instead of $j$, and hence
the $k$th terms agree in this case too. Finally, for
$k=j$, $G_j$ is new at $j$ in $\mes(\MM)$ by the
assumption, so $G_j$ is the smallest set in $\MM$
not containing $j$. Then $F$ is even smaller such
set, and so it comes to the $j$th position
of $\mes(\NN)$.
This concludes the proof of Proposition~\ref{PropCol}.
\proofend

\section{Collapsibility Versus Leray Number}

In this section we prove Theorem~\ref{t:}(b).

First we recall the notion of join of two simplicial
complexes $\K$ and $\L$. First, assuming that the vertex
sets $V(\K)$ and $V(\L)$ are disjoint, the \emph{join}
is the simplicial complex $\K*\L:=\{\sigma\cup\tau:
\sigma\in\K,\tau\in\L\}$ on the vertex set $V(\K)\cup V(\L)$.
If the vertex sets are not disjoint (as will be the case
in our application below), we first take isomorphic
copies of $\K$ and $\L$ with disjoint vertex sets and
then we form the join as above.

The next lemma shows that the Leray number behaves
nicely with respect to joins.

\begin{lemma}\label{joins-ler}
For every two nonempty simplicial complexes $\K$ and $\L$
we have $\ldim(\K*\L)=\ldim(\K)+\ldim(\L)$.
\end{lemma}

\begin{proof} 
This is a simple consequence of the K\"unneth formula for joins
\begin{equation}\label{e:kunneth}
\HrFree k{\scxf X * \scxf Y} = \bigoplus_{i+j = k-1}
\HrFree i{\scxf X} \otimes \HrFree j{\scxf Y}
\end{equation}
for any two simplicial complexes $\scxf X$ and $\scxf Y$,
where $\HrFree k .$ denotes the $k$-dimensional reduced
homology group (over $\kve$). The K\"unneth formula 
in this form can easily be
derived from \cite{Mun-eoat}, Example~4 (p.~349)
 and Exercise~3 (p.~373). 

For notational convenience we assume
$V(\K)\cap V(\L)=\emptyset$, and let 
$\ldim(\K)=k$ and $\ldim(\L)=\ell$. Then
there exist $A\subseteq V(\K)$
and $B\subseteq V(\L)$ such that $\HrFree {k-1}{\K[A]}\ne
0\ne \HrFree {\ell-1}{\L[B]}$. Since
$(\K*\L)[A\cup B]= \K[A]*\L[B]$, (\ref{e:kunneth})
shows that $\HrFree{k+\ell-1}{(\K*\L)[A\cup B]}\ne 0$,
and thus $\ldim(\K*\L)\ge k+\ell$.

On the other hand, every induced subcomplex of $\K*\L$
has the form $\K[A]*\L[B]$ for some $A\subseteq V(\K)$
and $B\subseteq V(\L)$, and if $\HrFree i{ \K[A]}=0$
for all $i\ge k$ and $\HrFree j{ \L[B]}=0$
for all $j\ge\ell$, (\ref{e:kunneth}) gives
$\HrFree s{\K[A]*\L[B]}=0$ for all $s\ge k+\ell$,
thus showing $\ldim(\K*\L)\le k+\ell$.
\end{proof}

\medskip

We cannot say how $\cdim(.)$ behaves under joins,
but we can do so for the following related quantity:
$$
\cols{\K} :=
\min\{d: \K \mbox{ has a $d$-collapsible 
    face}\}.
$$

\begin{lemma} \label{joins-coll} For every two
simplicial complexes $\K,\L$ we have
$\cols{\K * \L} = \cols{\K} + \cols{\L}$.
\end{lemma} 

\begin{proof} Again we assume
$V(\K)\cap V(\L)=\emptyset$.
It is easily checked that if $\sigma$ is a $k$-collapsible
face of $\K$ and $\tau$ is an $\ell$-collapsible
face of $\L$, then $\sigma\cup\tau$ is
a $(k+\ell)$-collapsible
face of $\K*\L$, which shows 
$\cols{\K*\L} \le \cols{\K} + \cols{\L}$.
On the other hand, \emph{every} $d$-collapsible
face of $\K*\L$ is of the form $\sigma\cup\tau$,
$\sigma\in \K$, $\tau\in\L$, and one can check
that  $\sigma$ is $k$-collapsible
and $\tau$ is $\ell$-collapsible for 
some $k,\ell$ with $k+\ell=d$. This gives the
reverse inequality.
\end{proof}

\heading{Proof of Theorem~\ref{t:}(b). }
We let $\K_0$ be the triangulation of the dunce hat
in Fig.~\ref{FigDH}. We have
$\ldim(\K_0)=2<\cols{\K_0}$ according to Wegner
\cite{w-dcnfc-75}, and actually $\cols{\K_0}=3$
because $\dim\K_0=2$. Then the join $\K$ of $d$ copies
of $\K_0$ satisfies $\ldim(K)=2d$ and $\cdim(\K)\ge
\cols{\K}=3d$ by Lemmas~\ref{joins-ler} and~\ref{joins-coll}.
\proofend


\bibliographystyle{alpha}
\bibliography{../../bib/cg.bib,../../bib/geom.bib}

\begin{thebibliography}{AKMM02}

\bibitem[AK92]{ak-pcs-92a}
N.~Alon and D.~Kleitman.
\newblock Piercing convex sets and the {Hadwiger Debrunner} $(p,q)$-problem.
\newblock {\em Adv. Math.}, 96(1):103--112, 1992.

\bibitem[AKMM02]{AKMM01}
N.~Alon, G.~Kalai, J.~Matou\v{s}ek, and R.~Meshulam.
\newblock Transversal numbers for hypergraphs arising in geometry.
\newblock {\em Adv. Appl. Math.}, 130:2509--2514, 2002.

\bibitem[Ame96]{a-spiht-96}
N.~Amenta.
\newblock A short proof of an interesting {Helly}-type theorem.
\newblock {\em Discrete Comput. Geom.}, 15:423--427, 1996.

\bibitem[B{\'a}r82]{b-gct-82}
I.~B{\'a}r{\'a}ny.
\newblock A generalization of {C}arath\'eodory's theorem.
\newblock {\em Discrete Math.}, 40:141--152, 1982.

\bibitem[DGK63]{DanzerGrunbaumKlee}
L.~Danzer, B.~Gr{\"u}nbaum, and V.~Klee.
\newblock Helly's theorem and its relatives.
\newblock In {\em Convexity}, volume~7 of {\em Proc. Symp. Pure Math.}, pages
  101--180. American Mathematical Society, Providence, 1963.

\bibitem[Eck93]{e-hrctt-93}
J.~Eckhoff.
\newblock Helly, {R}adon and {C}arath\'eodory type theorems.
\newblock In P.~M. Gruber and J.~M. Wills, editors, {\em Handbook of Convex
  Geometry}. North-Holland, Amsterdam, 1993.

\bibitem[Flo34]{f-undkd-}
A.~Flores.
\newblock {\"U}ber $n$-dimensionale {Komplexe} die im {$R_{2n+1}$} absolut
  selbstverschlungen sind.
\newblock {\em Ergeb. Math. Kolloq.}, 6:4--7, 1932/1934.

\bibitem[Hat01]{Hatcher}
A.~Hatcher.
\newblock {\em Algebraic Topology}.
\newblock Cambridge University Press, Cambridge, 2001.
\newblock Electronic version available at {\tt http://math.cornell.edu/\~{
  }hatcher\#AT1}.

\bibitem[Hel23]{h-umkkm-23}
E.~Helly.
\newblock {\"U}ber {Mengen} konvexer {K{\"o}rper} mit gemeinschaftlichen
  {Punkten}.
\newblock {\em Jahresbericht Deutsch. Math. Verein.}, 32:175--176, 1923.

\bibitem[Hel30]{h-usvam-30}
E.~Helly.
\newblock {\"U}ber {Systeme} von abgeschlossenen {Mengen} mit
  gemeinschaftlichen {Punkten}.
\newblock {\em Monaths. Math. und Physik}, 37:281--302, 1930.

\bibitem[Kal84]{Kalai-nececkhoff}
G.~Kalai.
\newblock Characterization of f-vectors of families of convex sets in {$R^d$}.
  {I}: Necessity of {E}ckhoff's conditions.
\newblock {\em Isr. J. Math.}, 48:175--195, 1984.

\bibitem[Kal86]{Kalai-suffeckhoff}
G.~Kalai.
\newblock {Characterization of $f$-vectors of families of convex sets in $R^d$.
  II: Sufficiency of Eckhoff's conditions}.
\newblock {\em J. Combin. Theory, Ser. A}, 41:167--188, 1986.

\bibitem[KL79]{kl-pgr-79}
M.~Katchalski and A.~Liu.
\newblock A problem of geometry in ${R}^n$.
\newblock {\em Proc. Amer. Math. Soc.}, 75:284--288, 1979.

\bibitem[KM05]{KalMes1}
G.~Kalai and R.~Meshulam.
\newblock {A topological colorful Helly theorem}.
\newblock {\em Adv. Math.}, 191(2):305--311, 2005.

\bibitem[KM07]{KalMes2}
G.~Kalai and R.~Meshulam.
\newblock Leray numbers of projections and a topological {H}elly type theorem.
\newblock Manuscript, The Hebrew University of Jerusalem, 2007.

\bibitem[Lov74]{Lovasz-colhelly}
L.~Lov{\'a}sz.
\newblock Problem 206.
\newblock {\em Matematikai Lapok}, 25:181, 1974.

\bibitem[Mat02]{Mat-dg}
J.~Matou\v{s}ek.
\newblock {\em Lectures on Discrete Geometry}.
\newblock Springer, New York, 2002.

\bibitem[Mun84]{Mun-eoat}
J. R. Munkres.
\newblock {\em Elements of Algebraic Topology}.
\newblock Addison-Wesley Pub., New York, 1984.


\bibitem[Rad21]{r-mkkde-21}
J.~Radon.
\newblock Mengen konvexer {K{\"o}rper}, die einen gemeinsamen {Punkt}
  enthalten.
\newblock {\em Math. Ann.}, 83:113--115, 1921.

\bibitem[vK32]{k-ker-32}
R.~E. van Kampen.
\newblock Komplexe in euklidischen {R\"a}umen.
\newblock {\em Abh. Math. Sem. Hamburg}, 9:72--78, 1932.
\newblock Berichtigung dazu, {\em ibid.} (1932) 152--153.

\bibitem[Weg75]{w-dcnfc-75}
G.~Wegner.
\newblock $d$-collapsing and nerves of families of convex sets.
\newblock {\em Arch. Math.}, 26:317--321, 1975.

\end{thebibliography}
\end{document}